\documentclass[a4paper,notitlepage,twoside,reqno,11pt]{amsart}

\usepackage{bbm,pifont,latexsym,amssymb}
\usepackage{anysize}
\marginsize{3.4cm}{3.4cm}{3cm}{3cm}

\usepackage{dcolumn,indentfirst}
\usepackage[hypertexnames=false,bookmarksopen=true,linktocpage=true,pdfstartview={XYZ null null 1.25}]{hyperref}
\usepackage{amsmath,amssymb,amscd,amsthm,amsfonts,mathrsfs}
\usepackage{color,graphicx,xcolor,graphics,subfigure,extarrows,caption2}
\usepackage{titletoc}

\usepackage{autonum}

\newtheorem{thm}{Theorem}[section]
\newtheorem{lem}[thm]{Lemma}
\newtheorem{prop}[thm]{Proposition}
\newtheorem{cor}[thm]{Corollary}

\newtheorem*{mainL}{Main Lemma}

\newtheorem{thmx}{Theorem}

\theoremstyle{definition}

\newcommand{\D}{\mathbb{D}}
\newcommand{\T}{\mathbb{T}}
\newcommand{\R}{\mathbb{R}}
\newcommand{\Z}{\mathbb{Z}}

\newcommand{\C}{\mathbb{C}}

\newcommand{\EC}{\widehat{\mathbb{C}}}

\newcommand{\MP}{\mathcal{P}}

\newcommand{\ME}{\mathcal{E}}

\newcommand{\diam}{\textup{diam}}
\newcommand{\dist}{\textup{dist}}

\newcommand{\Int}{\textup{int}}
\newcommand{\Log}{\textup{Log}}
\newcommand{\Ext}{\textup{ext}}
\newcommand{\rot}{\textup{rot}}

\newcommand{\re}{\textup{Re}\,}
\newcommand{\im}{\textup{Im}\,}

\makeatletter\@addtoreset{equation}{section}\makeatother

\begin{document}

\author[F. Yang]{Fei Yang}
\address{School of Mathematics, Nanjing University, Nanjing 210093, P. R. China}
\email{yangfei@nju.edu.cn}

\author[G. Zhang]{Gaofei Zhang}
\address{School of Mathematics, Nanjing University, Nanjing 210093, P. R. China}
\email{zhanggf@nju.edu.cn}

\author[Y. Zhang]{Yanhua Zhang}
\address{School of Mathematical Sciences, Qufu Normal University, Qufu 273165, P. R. China}
\email{zhangyh0714@qfnu.edu.cn}

\title[A new premodel for bounded type Siegel disks]{A new premodel for bounded type Siegel disks of some transcendental entire functions}

\begin{abstract}
In the 1980s, Douady and Sullivan conjectured that the Siegel disks of rational maps with degree at least two are Jordan domains.
Based on the study of Blaschke models, this conjecture was settled by Zhang under the condition that the Siegel disks are of bounded type, and by Zakeri for a large number of transcendental entire functions under the same condition.
The entire functions studied by Zakeri have exactly one finite asymptotic value, which is the center of the fixed Siegel disk. Ch\'{e}ritat, and later with Epstein, also studied the Siegel disk boundaries of such type transcendental entire functions by developing Ghys-like models.

In this paper, by introducing a new premodel, we prove that the fixed bounded type Siegel disk centered at the \textit{non-zero} point $\lambda-1$ of $f_\lambda(z)=e^{1-\lambda}ze^z$ is a quasi-disk whose boundary contains the critical point $-1$. The asymptotic value $0$ is not contained in the Siegel disk of these entire functions, which is essentially different from the previous results. Moreover, we prove that the Julia sets of these functions with a fixed bounded type Siegel disk are locally connected.
\end{abstract}



\date{\today}

\maketitle

\section{Introduction}

\subsection{Backgrounds}

In holomorphic dynamics, a \textit{Siegel disk} is a simply connected Fatou component in which the dynamics of a holomorphic map is conjugate to an irrational rotation of the unit disk.
More precisely, if a holomorphic map $f$ has a fixed point at the $z_0$ with multiplier $e^{2\pi i\theta}$ for an irrational number $\theta\in \mathbb{R}\setminus \mathbb{Q}$ and $f$ is locally linearizable at $z_0$, then the maximal domain in which the linearization exists is called the Siegel disk of $f$ centered at $z_0$.
It is well known that if $\theta$ is a Brjuno number, then $f$ has a Siegel disk centered at $z_0$ (\cite{Sie42}, \cite{Brj71}). In particular, this is true for \textit{bounded type} irrational number $\theta$, i.e., the continued fraction expansion $[a_0;a_1,a_2,\cdots,a_n,\cdots]$ of $\theta$ satisfies $\sup_{n\geq 1} \{a_n\}<\infty$.

The dynamics in Siegel disks is relatively easy to be understood. People are interested in the geometric and topological properties of their boundaries. In the early 1980s, Douady and Sullivan conjectured that the Siegel disk boundaries of rational maps with degree at least two are Jordan curves (\cite{Dou83}, \cite{Rog92a}).
In 2011, the second author of this paper settled this conjecture under the assumption that the rotation number $\theta$ is of bounded type \cite{Zha11}. Actually he proved more: the boundaries of these Siegel disks are quasi-circles passing through at least one critical points. This recovers the previous results obtained by Douady-Herman for quadratic polynomials \cite{Dou87}, \cite{Her87}, Zakeri for cubic polynomials \cite{Zak99}, and Shishikura for all polynomials \cite{Shi01}.
One may also refer to \cite{ABC04}, \cite{PZ04}, \cite{BC07}, \cite{Che25}, \cite{SY25} for related results on this conjecture when $\theta$ is not of bounded type.

The analogous question for transcendental entire functions is more subtle: Siegel disks in the transcendental setting could be unbounded, and their boundaries may exhibit more complicated topological behavior than those of polynomials or rational maps.
Hence for a transcendental entire function with a Siegel disk $\Delta$, a natural question is under what conditions $\Delta$ is bounded, and further, whether its boundary is a Jordan curve.

In 2010, building on Shishikura's pioneering method of quasiconformal surgery, Zakeri proved the following result.

\begin{thm}[\cite{Zak10}]\label{thm:Zakeri}
Let $f$ be a non-linear entire function of the form
\begin{equation}
f(z)=P(z)e^{Q(z)},
\end{equation}
where $P$ and $Q$ are polynomials of degrees $p\geq 1$ and $q\geq 0$ respectively. If $f$ has a bounded type Siegel disk centered at the origin, then the boundary of this Siegel disk is a quasi-circle containing at least one critical point of $f$.
\end{thm}

For polynomials (the case $q=0$), Theorem \ref{thm:Zakeri} recovers Shishikura's result \cite{Shi01}. It also extends Geyer's work on the map $z\mapsto e^{2\pi i\theta}ze^{z}$ \cite{Gey01} and Keen-Zhang's result for the maps of the form $z\mapsto (e^{2\pi i\theta}z+az^2)e^z$ \cite{KZ09}, where $\theta$ is of bounded type.
Zakeri's proof is based on the study of a family of analytic maps
\begin{equation}
G(z)=\lambda z^p B(z)\exp\big(\alpha(z)-\overline{\alpha(1/\overline{z})}\big):\C\setminus\{0\}\to\EC,
\end{equation}
where $|\lambda|=1$, $B$ is a Blaschke product of degree $p-1$ and $\alpha$ is a polynomial of degree $q$ with $\alpha(0)=0$.
These maps serve as models of the maps in Theorem \ref{thm:Zakeri}. The family of these analytic maps is compact (the same property of finite Blaschke products studied by Shishikura), which is essential since this guarantees that the corresponding critical circle homeomorphisms can be uniformly quasisymmetrically linearizable. This in turn implies that the invariant curves in the Siegel disk are uniform quasi-circles (i.e., the dilatation has a uniform upper bound), which implies that the boundary of the Siegel disk is a quasi-circle.

However, Zakeri emphasized in \cite{Zak10} that the result of Theorem \ref{thm:Zakeri} may not be true if the center of the Siegel disk of $f(z)=P(z)e^{Q(z)}$ is not at the origin.
Indeed, for the exponential map $z\mapsto \lambda e^{z-\lambda}$ with $\lambda=e^{2\pi i\theta}$, the fixed bounded type Siegel disk centered at $\lambda$ is unbounded \cite{Her85} and its boundary is not locally connected \cite{BD00b}. Note that the Siegel center $\lambda$ has infinitely many preimages.
If one applies the strategy of Zakeri (or Shishikura) to construct a family of Blaschke-like models, then these models have infinitely many zeros and the quasisymmetric constant of the linearizing functions would not be controlled uniformly.

We would like to mention that for the sine family $z\mapsto\lambda\sin z$, although the fixed bounded type Siegel disk $\Delta$ centered at $0$ has infinitely many preimages, the second author of this paper invented a periodic premodel and proved that $\partial\Delta$ is a quasi-circle passing through exactly two critical points $\frac{\pi}{2}$ and $-\frac{\pi}{2}$ \cite{Zha05}. See also \cite{Yan13} and \cite{ZFS20} for general situations.

We also note that Ch\'{e}ritat, and later with Epstein, proved that a family of finite type transcendental entire functions with at most two singular values in $\C$ have a Siegel disk whose boundary is a quasi-circle containing exactly one critical point by developing the Ghys-like model, and they also required that the center of the Siegel disk is a singular value \cite{Che06}, \cite{CE18}.

\subsection{Main results}

In this paper, we consider a family of non-periodic transcendental entire functions with a bounded type Siegel disk which does not contain any singular values and prove that the boundaries of these Siegel disks are quasi-circles passing through one critical point by introducing a new premodel.

\begin{thmx}\label{thm:main}
Let $\theta$ be a bounded type irrational number and $\lambda=e^{2\pi i\theta}$. Then
\begin{equation}
f_{\lambda}(z)=e^{1-\lambda} ze^z
\end{equation}
has a fixed Siegel disk centered at the non-zero point $\lambda-1$ with rotation number $\theta$ whose boundary is a quasi-circle containing the critical point $-1$.
\end{thmx}

\begin{figure}[!htpb]
  \setlength{\unitlength}{1mm}
  \centering
 \includegraphics[width=.95\textwidth]{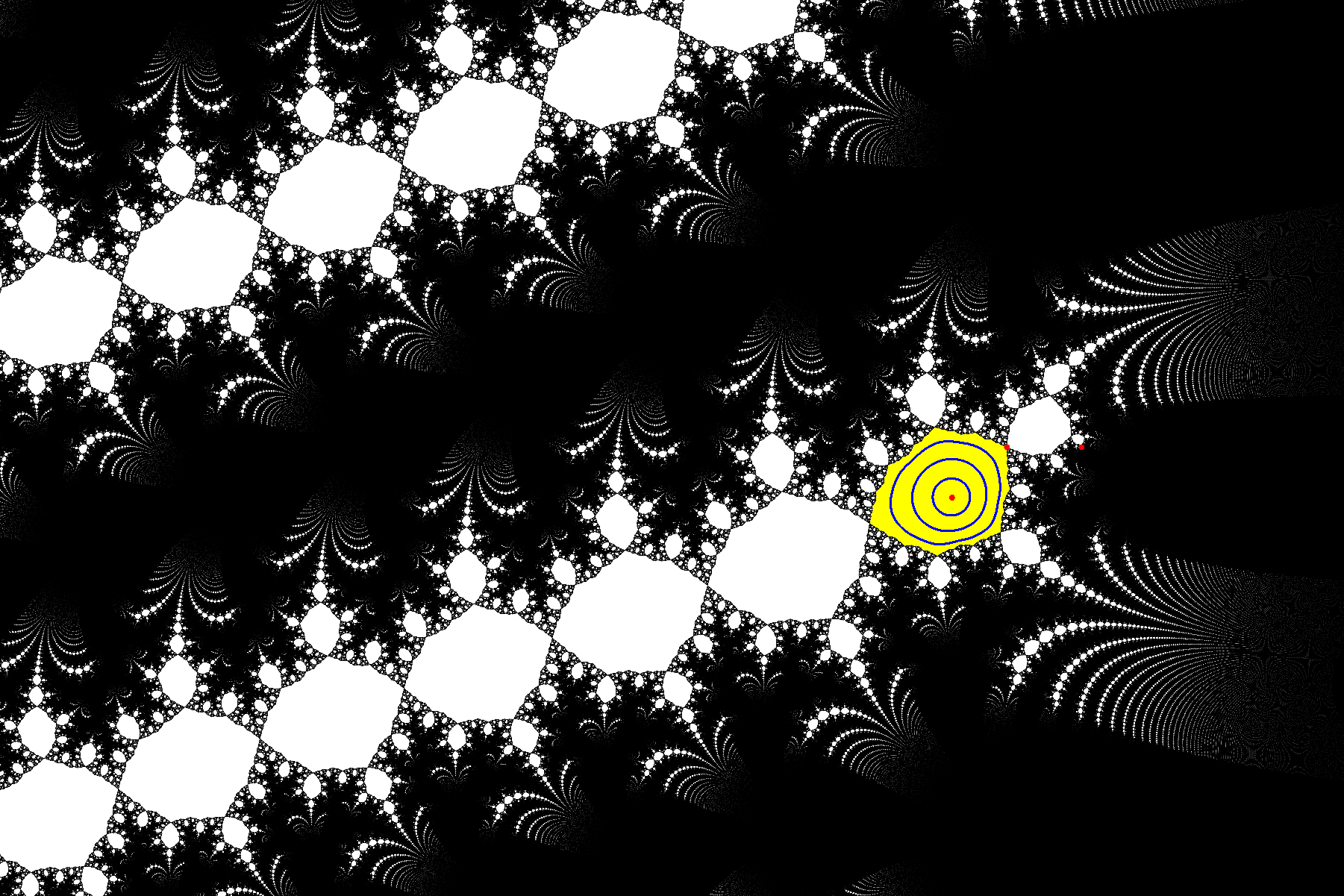}
 \put(-25,43){\textcolor{white}{$0$}}
 \put(-33.5,45.1){\small{$-1$}}
  \caption{The dynamical plane of $f_{\lambda}(z)=e^{1-\lambda}ze^{z}$, where $\lambda=e^{2\pi i\theta}$ with $\theta=(\sqrt{5}-1)/2$. The Siegel disk is colored yellow with invariant curves surrounding the center $\lambda-1$. The critical point $-1$ and the asymptotic value $0$ are also marked. All white regions are preimages of the Siegel disk and the Julia set is the black part. The bounary of the Siegel disk is a quasi-circle and the Julia set is locally connected.}
\label{Fig:Julia}
\end{figure}

The family of entire functions $g_a(z)=a ze^z$, where $a\in\C\setminus\{0\}$, was studied by Baker first, who used this family to produce Julia sets which are the whole complex plane \cite{Bak70}. In 1995, Fagella studied the parameter space of this family and obtained the connection between this family and Arnold standard family for some special $a$ \cite{Fag95}. Later, Geyer used the complexification of the Arnold standard family to serve as the model of $g_a$ with $a=e^{2\pi i\theta}$ and proved that if $\theta$ is of bounded type, then the Siegel disk of $g_a$ centered at the origin is bounded by a quasi-circle containing the critical point $-1$ \cite{Gey01}.
For $f_\lambda$ in Theorem \ref{thm:main}, the asymptotic value $0$ is a repelling fixed point, which is not in the Siegel disk. See Figure \ref{Fig:Julia}.

Recently, the local connectivity of the Julia sets of some transcendental entire functions with bounded type Siegel disks are proved by assuming the Siegel disks are bounded by quasi-circles containing at least one critical point. Based on Theorem \ref{thm:main}, as an application of the criterion in \cite{YZZ25}, we obtain the following result.

\begin{thmx}\label{thm:JLC}
Let $f_\lambda$ be the transcendental entire function with a bounded type Siegel disk in Theorem \ref{thm:main}. Then the Julia set of $f_\lambda$ is locally connected.
\end{thmx}

In \cite[Theorem 1.3]{YZZ25}, we prove the existence of $|a|>1$ such that $g_a(z)=a z e^z$ is polynomial-like renormalizable, has a $2$-cycle of bounded-type Siegel disks whose boundary components are quasi-circles, and has a locally connected Julia set. See \cite{BF10} and \cite{KN22b} for more examples of transcendental entire functions which are polynomial-like renormalizable and have quasiconformal copies of filled Julia sets of polynomials with bounded type Siegel disks.
The functions in Theorem \ref{thm:main} are not polynomial-like renormalizable (Corollary \ref{cor:non-renorm}).

\subsection{Idea of the proof}

\begin{figure}[!htpb]
  \setlength{\unitlength}{1mm}
  \centering
  \includegraphics[width=0.85\textwidth]{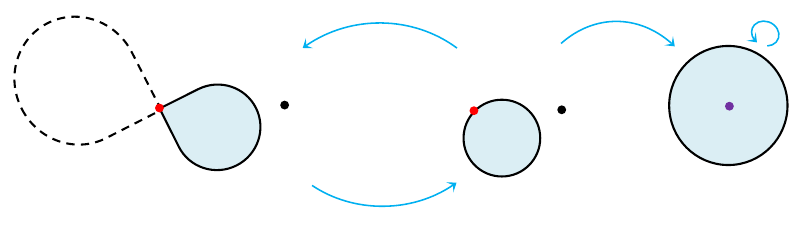}
  \put(-100.5,12.5){$-1$}
  \put(-89.5,14){$U$}
  \put(-79,13.5){$0$}
  \put(-56,18){$-\frac{1}{e}$}
  \put(-37.2,13){$0$}
  \put(-46.5,12){$V$}
  \put(-40.5,5){$\Gamma$}
  \put(-12,13){$0$}
  \put(-15,21){$\D$}
  \put(-64,32.5){$\phi$}
  \put(-65,0){$g$}
  \put(-29,32){$\Psi$}
  \put(-6,32){$R_\theta$}
  \caption{The construction of the premodel. The key point is locating the position of the round circle $\Gamma$.}
\label{Fig:model}
\end{figure}

Let us give a quick outline of the proof of Theorem \ref{thm:main}. Our proof is inspired by Ch\'{e}ritat's construction of Ghys-like premodel \cite{Che06} and we will mention the main difference in the construction.
The map
\begin{equation}
g(z):=g_1(z)= ze^z
\end{equation}
has exactly two singular values $0$ and $-\frac{1}{e}$ in $\C$. Let $\Gamma$ be a round circle in $\C\setminus\{0\}$ passing through the critical value $-\frac{1}{e}$ such that the asymptotic value $0$ lies in the exterior of $\Gamma$. We use $V$ to denote the Jordan domain bounded by $\Gamma$ and let $U$ be a bounded component of $g^{-1}(V)$ whose boundary contains the critical point $-1$. Note that $g:\partial U\to\Gamma$ is a homeomorphism. For each given bounded type $\theta$, let $\phi:\EC\setminus\overline{V}\to\EC\setminus\overline{U}$ be the unique conformal isomorphism fixing $\infty$ such that
\begin{equation}
g\circ \phi|_{\Gamma}: \Gamma\rightarrow \Gamma
\end{equation}
has rotation number $\theta$. Since $g\circ \phi|_\Gamma$ can be extended analytically to a neighborhood of $\Gamma$, it follows that $g\circ \phi|_\Gamma$ is quasisymmetrically conjugate to the rigid rotation $R_{\theta}(\zeta)=e^{2\pi i\theta}\zeta$ on the unit circle. We extend the conjugacy to a homeomorphism $\Psi:\overline{V}\to \overline{\D}$ such that  $\Psi:V\to \D$ is quasiconformal. Let $F:\C\to\C$ be the quasi-regular map defined by $g\circ \phi$ in $\C\setminus V$ and $\Psi^{-1}\circ R_\theta\circ \Psi$ in $V$. Then one can find a quasiconformal mapping $\Phi:\EC\to\EC$ fixing $0$ and $\infty$ such that $f:=\Phi\circ F\circ\Phi^{-1}$ is transcendental entire function having a Siegel disk whose boundary is a quasi-circle containing a critical point. See Figure \ref{Fig:model}.

To prove that $f$ has the form $f_\lambda(z)=e^{1-\lambda}ze^{z}$ in Theorem \ref{thm:main}, it suffices to prove $\phi(0)=0$ (see \S\ref{subsec:pf-A}). However, the conformal map $\phi:\EC\setminus\overline{V}\to\EC\setminus\overline{U}$ has been uniquely determined by the requirement that $g\circ \phi|_{\Gamma}: \Gamma\rightarrow \Gamma$ has rotation number $\theta$. To overcome this difficulty, we consider a family of round circles $\{\Gamma_t\}_t$ and prove that there exists a suitable round circle $\Gamma=\Gamma_{t_0}$ in $\C$ such that the corresponding $\phi=\phi_{t_0}$ satisfies $\phi(0)=0$. In particular, these two conditions are compatible (see the \hyperref[lem:main]{Main Lemma} in \S\ref{sec:proof-main} and its proof in \S\ref{sec:model}). This is the main difference between our construction and Cheritat's. In \cite{Che06}, the curve $\Gamma$ is chosen as the round circle centered at $0$ with radius $\frac{1}{e}$ and the property $\phi(0)=0$ is not needed there.

\medskip
\noindent\textbf{Organization of the paper.}
In \S\ref{sec:proof-main} we use quasiconformal surgery to prove Theorem \ref{thm:main} by assuming the \hyperref[lem:main]{Main Lemma} (i.e., the existence of the circle $\Gamma$). We also prove Theorem \ref{thm:JLC} based on Theorem \ref{thm:main} and a weak expansion property near the bounded type Siegel disks of transcendental entire functions obtained in \cite{YZZ25}. In \S\ref{sec:model}, we prove the existence of $\Gamma$ by deforming the round circles $\{\Gamma_t\}_t$.

\medskip
\noindent\textbf{Notations.}
We use $\mathbb{C}$ and $\widehat{\mathbb{C}}=\mathbb{C}\cup\{\infty\}$ to denote the complex plane and the Riemann sphere respectively. Let $\D(a,r)=\{z\in\C:|z-a|<r\}$ be the Euclidean disk centered at $a\in\C$ with $r>0$. For short, we denote $\D_r:=\D(0,r)$ and $\T_r:=\partial\D_r$. In particular, $\D:=\D_1$ is the unit disk and $\T:=\partial\D$ is the unit circle.
For a Jordan curve $\gamma$ in $\C$, we use $\gamma^{\Ext}$ to denote the component of $\EC\setminus\gamma$ containing $\infty$ and use $\gamma^{\Int}$ to denote the other.
Let $X$ be $\C$, $\EC$ or a hyperbolic domain in $\C$. We use $\dist_X(\cdot,\cdot)$ and $\diam_X(\cdot)$, respectively, to denote the distance and diameter with respect to the Euclidean, spherical or hyperbolic metric.

\medskip
\noindent\textbf{Acknowledgments.}
This work was supported by NSFC (Grant No. 12571093) and NSF of Shandong Province (Grant No. ZR2023MA044).

\section{Proofs of the main results}\label{sec:proof-main}

In this section, we give the detailed construction of a premodel and prove Theorems \ref{thm:main} and \ref{thm:JLC} by assuming the existence of a round circle with nice properties.

\subsection{Some basic properties}

Let $U_0$ and $U_1$ be two Jordan domains in $\C$ such that $U_0$ is compactly contained in $U_1$. The map $f:U_0\to U_1$ is called a \textit{polynomial-like map} of degree $d\geq 2$ if $f$ is a proper holomorphic surjection of degree $d$. Let $f$ be a transcendental entire function. If there exists $p\geq 1$ such that $f^{\circ p}:U_0\to U_p$ is a polynomial-like map of degree at least two, then $f$ is called \textit{$p$-polynomial-like renormalizable} (see \cite{DH85b}, \cite{McM94b}).
We first prove the following result.

\begin{lem}\label{lem:non-renorm}
Let $g_a(z)=a ze^z$ with $a\in\C\setminus\{0\}$. Then $g_a$ is not $1$-polynomial-like renormalizable.
\end{lem}

\begin{proof}
Assume that $g_a$ is $1$-polynomial-like renormalizable. Then there exist two Jordan disks $U_0$ and $U_1$ such that $U_0$ is compactly contained in $U_1$ and $g_a:U_0\to U_1$ is a quadratic-like map (i.e., a polynomial-like map of degree two) since $-1$ is the only critical point and the local degree of $-1$ is two.
We claim that $0\not\in U_1$. Indeed, otherwise $g_a:U_0\to U_1$ is proper of degree one since $g_a^{-1}(0)=\{0\}$ and the local degree of $g_a$ at $0$ is one, which is a contradiction. Hence $0\not\in U_1$.

Let $\Log(\cdot)$ be a branch of the inverse of the exponential map $z\mapsto e^z$ defined in $U_1$. Then the map $G(z):=\Log(g_a(z))-\Log(z)-z$ is well-defined and holomorphic in $U_0$. Note that $e^{G(z)}=a$ is a constant. Hence $G$ is also a constant and there exists $b\in\C$ satisfying $e^b=a$ such that
\begin{equation}
\Log(g_a(z))-\Log(z)=z+b, \text{ for all } z\in U_0.
\end{equation}
If $z_0$ is a fixed point of $g_a$ in $U_0\subset\C\setminus\{0\}$, then $z_0=-b\neq 0$ and hence different fixed points of $g_a$ cannot lie in $U_0$ at the same time.
However, $g_a:U_0\to U_1$ is a quadratic-like map. Hence $g_a$ has exactly one fixed point in $U_0$ with multiplicity two. This implies that $z_0$ is a parabolic fixed point with multiplier $1$. A direct calculation shows that $g_a(z)=ze^z$ and $z_0=0$, which is a contradiction. The proof is complete.
\end{proof}

\begin{cor}\label{cor:non-renorm}
The map $f_\lambda(z)=e^{1-\lambda} ze^z$ with $|\lambda|\leq 1$, is not polynomial-like renormalizable.
\end{cor}

\begin{proof}
Assume that $f_\lambda$ is $p$-polynomial-like renormalizable, where $p\geq 1$. Note that $f_\lambda$ has exactly one critical point $-1$ and one non-repelling fixed point at $\lambda-1$. There exist Jordan domains $U_0$ and $U_p$ such that $f^{\circ p}:U_0\to U_p$ is a quadratic-like map and $-1$, $\lambda-1\in U_0$. Since $\lambda-1$ is a fixed point, we have $p=1$. By Lemma \ref{lem:non-renorm}, this is impossible. Hence $f_\lambda$ is not polynomial-like renormalizable.
\end{proof}

For integers $p,q\geq 0$, let $\mathcal{E}^{p,q}$ be the family of nonlinear entire functions of the form
\begin{equation}\label{equ:f-P-Q}
f(z)=P(z)e^{Q(z)},
\end{equation}
where $P$ and $Q$ are polynomials of degrees $p$ and $q$ respectively.
To make a definite choice, for each $f\in\mathcal{E}^{p,q}$ we use the unique representation of the form \eqref{equ:f-P-Q} in which $Q(0) = 0$.
We first recall some basic properties of the maps in the family $\ME^{p,q}$ (see \cite[\S 3]{Zak10}).
The \textit{growth order} of an entire map $f: \mathbb{C} \rightarrow \mathbb{C}$ is defined by
\begin{equation}\label{equ:growth-order}
\limsup _{r \rightarrow+\infty} \frac{\log \log M(f, r)}{\log r},
\end{equation}
where $M(f, r):=\sup _{|z|=r}|f(z)|$. In particular, the growth order of every map in $\mathcal{E}^{p, q}$ is $q$ and we have $g_a(z)=aze^z\in\ME^{1,1}$ for all $a\in\C\setminus\{0\}$.

Note that the growth order in \eqref{equ:growth-order} can be also defined for quasi-regular maps defined from $\C$ to itself.
The following result is a slight modification of \cite[Corollary 3.2]{Zak10} and the proof is completely the same (by applying the H\"{o}lder continuity of quasiconformal mappings).

\begin{lem}\label{lem:E-rigidity}
Suppose $f:\C\to\C$ is an entire function and there are quasiconformal mappings $\varphi,\hat{\varphi}:\C\to\C$ which fix zero such that $\varphi^{-1}\circ f\circ\hat{\varphi}$ has finite growth order, $p$ zeros and $p+q-1$ critical points counting multiplicities. Then $f\in\ME^{p,q}$.
\end{lem}

We need the following result to study the preimages of simply connected domains (see also \cite[\S 3]{Zak10} for partial results on $f\in\ME^{p,q}$).

\begin{prop}[{\cite[Proposition 2.8]{BFR15} and \cite[Lemma 3.1]{HXY25}}]\label{prop:E-preimage}
Let $f$ be a transcendental entire function, $V \subset \mathbb{C}$ be a simply connected domain, and $\widetilde{V}$ be a component of $f^{-1}(V)$. Then either
\begin{itemize}
\item[(1)] $f: \widetilde{V} \rightarrow V$ is a proper map and hence has finite degree; or
\item[(2)] $f^{-1}(w) \cap \widetilde{V}$ is infinite for every $w \in V$, with at most one exception. In this case, either $\widetilde{V}$ contains an asymptotic curve corresponding to an asymptotic value in $V$, or $\widetilde{V}$ contains infinitely many critical points.
\end{itemize}
In case (1), if $V$ is a bounded Jordan domain such that $\overline{V}$ contains no asymptotic value and at most one critical value, then $\widetilde{V}$ is also a bounded Jordan domain.
\end{prop}

As an immediate corollary of Proposition \ref{prop:E-preimage}, we have the following result.

\begin{cor}\label{cor:preimage}
Let $g(z)=ze^z$ and $V$ be a domain in $\C\setminus\{0\}$ with $0\not\in\overline{V}$ whose boundary $\partial V$ is a real-analytic Jordan curve containing $-1/e$. Then $g^{-1}(V)$ has a component $U$ whose boundary is a Jordan curve containing the critical point $-1$ and is real-analytic except at $-1$. Moreover, $g:U\to V$ is conformal.
\end{cor}

\subsection{Construction of the premodel}

In the rest of this subsection, we fix the bounded type irrational number $\theta$.
Let $\Gamma$ be a round circle in $\C\setminus\{0\}$ passing through the critical value $-\frac{1}{e}$ of $g(z)=ze^z$ such that the asymptotic value $0$ lies in the exterior of $\Gamma$. Let $V$ be the bounded Jordan domain bounded by $\Gamma$.
By Corollary \ref{cor:preimage}, $g:\overline{U}\to \overline{V}$ is a homeomorphism, where $U$ is a bounded component of $g^{-1}(V)$ whose boundary is a Jordan curve containing the critical point $-1$ and is real-analytic except at $-1$. Let
\begin{equation}
\phi:\EC\setminus\overline{V}\to\EC\setminus\overline{U}
\end{equation}
be a conformal map fixing $\infty$. Note that $\phi$ can be extended to homeomorphism $\phi:\EC\setminus V\to\EC\setminus U$. Then $g\circ\phi|_{\Gamma}:\Gamma\to\Gamma$ is also an orientation-preserving homeomorphism. See Figure \ref{Fig:model}.

The most crucial ingredient in the construction of the premodel for the maps in Theorem \ref{thm:main} is the following result, whose proof will be given in \S\ref{sec:model}.

\begin{mainL}\label{lem:main}
There exists a round circle $\Gamma\subset\C\setminus\{0\}$ such that the rotation number of $g\circ\phi|_{\Gamma}:\Gamma\to\Gamma$ is $\theta$, and the conformal map $\phi:\EC\setminus\overline{V}\to\EC\setminus\overline{U}$ fixes $\infty$ and $0$.
\end{mainL}

In the rest of this subsection, we fix the choice of the circle $\Gamma$ and the conformal map $\phi$ in the \hyperref[lem:main]{Main Lemma}.
By the Schwarz reflection principle, we have the following result.

\begin{lem}\label{lem:analytic-extension}
The circle homeomorphism $g\circ \phi|_\Gamma:\Gamma\to\Gamma$ can be extended to an analytic map defined in $\C\setminus\{0\}$, and this map has local degree three at $\phi^{-1}(-1)$, which is the unique critical point on $\Gamma$.
\end{lem}

The following \textit{Herman-\'{S}wi\c{a}tek Theorem} gives a necessary and sufficient condition for an orientation-preserving homeomorphism of the unit circle to be quasisymmetrically conjugate to the rigid rotation (see \cite{Her87}, \cite{Swi98}, \cite{Pet00}).

\begin{thm}[{Herman-\'{S}wi\c{a}tek}]\label{thm:Herman-Swiatek}
Let $f:\T\to\T$ be a real-analytic critical circle homeomorphism of rotation number $\alpha$. Then $f$ is quasisymmetrically conjugate to $R_{\theta}(\zeta)=e^{2\pi i\alpha}\zeta$  if and only if $\alpha$ is of bounded type.
\end{thm}

Note that the result of Lemma \ref{thm:Herman-Swiatek} also holds if the unit circle $\T$ is replaced by any other round circle in $\C$.
By the \hyperref[lem:main]{Main Lemma}, Lemma \ref{lem:analytic-extension} and Theorem \ref{thm:Herman-Swiatek}, $g\circ \phi|_\Gamma:\Gamma\to\Gamma$ is a real-analytic critical circle homeomorphism with bounded type rotation number $\theta$, which is quasisymmetrically conjugate to $R_{\theta}$. Thus, there exists a unique quasisymmetric
homeomorphism $\psi:\Gamma\to\T$  such that
\begin{equation}\label{equ:q0}
\psi\circ g\circ \phi\circ\psi^{-1}=R_\theta \text{\quad and\quad} \psi(q_0)=-1, \text{ where }q_0=\phi^{-1}(-1)\in\Gamma.
\end{equation}
By the Beurling-Ahlfors extension, there exists a homeomorphism $\Psi:\overline{V}\to\overline{\mathbb{D}}$ which is quasiconformal in $V$ such that $\Psi|_{\Gamma}=\psi$.
Set
\begin{equation}\label{equ:premodel}
F(z):=
\begin{cases}
   g \circ \phi(z), & z \in\C\setminus V, \\
  \Psi^{-1} \circ R_{\theta} \circ\Psi(z), & z \in V.
\end{cases}
\end{equation}
Then $F:\C\to\C$ is quasi-regular, which is the premodel of the maps in Theorem \ref{thm:main}.

\subsection{Proof of Theorem \ref{thm:main}}\label{subsec:pf-A}

In this subsection we give a proof of Theorem \ref{thm:main} by assuming the \hyperref[lem:main]{Main Lemma}.

\begin{proof}[Proof of Theorem \ref{thm:main}]
We define an $F$-invariant complex structure $\mu$ as follows. Let $\mu_0$ be the standard complex structure of the plane represented by the zero Beltrami
differential.
For $z\in V$, let $\mu(z)$ be the complex structure given by $\Psi^*\mu_0$.
Note that $F:V\to V$ is quasiconformal and $R_\theta^*\mu_0=\mu_0$. In the disk $V$ we have
\begin{equation}
F^*\mu=( \Psi^{-1} \circ R_{\theta} \circ\Psi)^*\mu=\Psi^*R_\theta^*(\Psi^{-1})^*\mu=\Psi^*R_\theta^*\mu_0=\Psi^*\mu_0=\mu.
\end{equation}
This means that $\mu$ is $F$-invariant in $V$.
For $z\in \C\setminus V $, the forward orbit of $z$ under $F$ either enters $V$ or does not.  In the former case, there exists a smallest integer $m\geq 1$ such that $w=F^{\circ m}(z)\in V$ and let $\mu$ be defined as the pullback of the complex structure at $w$ by $F^{\circ m}$.
In the latter case, we define $\mu(z)=\mu_0$. Clearly, $\mu$ is $F$-invariant and bounded (i.e., $\|\mu\|_\infty<1$).

By the measurable Riemann mapping theorem, there exists a unique quasiconformal homeomorphism $\Phi: \widehat{\mathbb{C}} \to \widehat{\mathbb{C}}$ integrating the complex structure $\mu$ such that $\Phi^*\mu_0=\mu$ and moreover,
\begin{equation}\label{equ:norm-Phi}
\Phi(\infty)=\infty,\quad\Phi(0)=0 \text{\quad and\quad}\Phi(q_0)=-1,
\end{equation}
where $q_0=\phi^{-1}(-1)\in\Gamma$ (see \eqref{equ:q0}).
Therefore, $f:=\Phi \circ F \circ \Phi^{-1}$ is an entire function.

By \eqref{equ:premodel}, $F$ has finite growth order since $g$ does. Moreover, $F$ has exactly one zero and one critical point. Hence by Lemma \ref{lem:E-rigidity}, we have $f\in\ME^{1,1}$.
By the \hyperref[lem:main]{Main Lemma}, we have $F(0)=0$ and hence $f(0)=0$. Thus $f$ has the form
\begin{equation}
f(z)=a ze^{bz} \text{\quad where } a,b\in\C\setminus\{0\}.
\end{equation}
Note that $f$ has exactly one critical point $-1$ by \eqref{equ:norm-Phi}. Hence $f(z)=a ze^z$.

By the construction, $f$ has a Siegel disk centered at a non-zero fixed point $z_0$ and its multiplier is $\lambda=e^{2\pi i\theta}$. A direct calculation shows that $f(z)=f_\lambda(z)=e^{1-\lambda}ze^z$. In particular, the boundary of the Siegel disk of $f_{\lambda}$ is a quasi-circle passing through the critical point $-1$. The proof of Theorem \ref{thm:main} is complete.
\end{proof}

\subsection{Proof of Theorem \ref{thm:JLC}}\label{subsec:LC}

In this subsection, we prove the Julia sets of the maps in Theorem \ref{thm:main} are locally connected, by using a weak expansion property near the bounded type Siegel disks of transcendental entire functions obtained in \cite{YZZ25}.
For a survey of results on the local connectivity of Julia sets of transcendental entire functions, see \cite{YZZ25}.

For a transcendental entire function $f:\C\to\C$, let $S(f)$ be the set of \textit{singular values} of $f$, i.e., the closure of critical values and asymptotic values of $f$ in $\C$. The \textit{postsingular set} of $f$ is
\begin{equation}\label{equ:P-f}
\MP(f):=\overline{\bigcup_{n\geq 0}f^{\circ n}\big(S(f)\big)}.
\end{equation}
We first determine all possible Fatou components of the maps in Theorem \ref{thm:main}.

\begin{lem}\label{lem:Fatou-finite}
Let $f_\lambda(z)=e^{1-\lambda} ze^z$, where $\lambda=e^{2\pi i\theta}$ and $\theta$ is of bounded type. Then the Fatou set $F(f_\lambda)$ is the union of the bounded type Siegel disk $\Delta_\lambda$ centered at $\lambda-1$ and its iterated preimages.
Each Fatou component of $f_\lambda$ is a bounded Jordan domain. Moreover, every periodic cycle in $J(f_\lambda)$ is repelling.
\end{lem}

\begin{proof}
Since $f_\lambda$ has only two singular values in $\C$, we know that $f_\lambda$ has no wandering domains and Baker domains (see \cite{GK86}, \cite{EL92}, \cite{Ye26}).
Note that the asymptotic value $0$ is a repelling fixed point and the critical value $-\frac{1}{e}=f_\lambda(-1)$ is contained in the boundary of the Siegel disk $\Delta_\lambda$ centered at $\lambda-1$. This implies that $f_\lambda$ has neither attracting nor parabolic basins.
By \cite[Theorem 7, p.\,164]{Ber93}, every Siegel disk boundary and Cremer point is accumulated by different points of the postsingular set $\MP(f_\lambda)=\partial\Delta_\lambda\cup\{0\}$. Therefore, every periodic cycle in $J(f_\lambda)$ is repelling and $f_\lambda$ has exactly one Siegel disk $\Delta_\lambda$.
Thus $F(f_\lambda)$ is the union of the bounded type Siegel disk $\Delta_\lambda$ centered at $\lambda-1$ and its iterated preimages.
By Theorem \ref{thm:main}, $\partial\Delta_\lambda$ is a quasi-circle. Note the Fatou set of $f_\lambda$ contains no critical point and the restriction of $f_\lambda$ on each Fatou component is conformal. By Proposition \ref{prop:E-preimage}, all Fatou components of $f_\lambda$ are bounded Jordan domains. The proof is complete.
\end{proof}

It is known that a continuum cannot fail to be locally connected only at a single point. Therefore, the Julia set $J(f)$ of a transcendental entire function $f$ is locally connected in $\C$ if and only if $J(f)\cup\{\infty\}$ is locally connected in $\EC$.
To prove the local connectivity of a Julia set in $\EC$, the following criterion is useful (see {\cite[Theorem 4.4, p.\,113]{Why42}}).

\begin{lem}[LC criterion]\label{lem:LC-criterion}
A connected compact subset $X$ in $\EC$ is locally connected if and only if the following two conditions hold:
\begin{enumerate}
\item  the boundary of every component of $\EC\setminus X$ is locally connected; and
\item  for any given $\varepsilon>0$, there are only finitely many components of $\EC\setminus X$ whose spherical diameters are greater than $\varepsilon$.
\end{enumerate}
\end{lem}

We will use the following result to prove the local connectivity of the Julia sets of the maps in Theorem \ref{thm:main}.

\begin{thm}[{\cite[Theorem 2.1]{YZZ25}}]\label{thm:main-pre}
Let $f$ be a transcendental entire function having a fixed bounded type Siegel disk $\Delta$ with quasi-circle boundary in $\C$. Suppose
\begin{enumerate}
  \item $\dist_{\C}(\MP(f)\setminus\partial{\Delta},\partial{\Delta})>0$;
  \item there are only finitely many critical values and no asymptotic values whose forward orbits intersect $\partial\Delta$; and
  \item the local degrees of the critical points whose forward orbits intersect $\partial\Delta$ have a uniform upper bound.
\end{enumerate}
Then for any $\varepsilon>0$ and bounded Jordan domain $V_0\subset \C\setminus\overline{\Delta}$ with $\overline{V}_0 \cap \MP(f)\neq\emptyset$ and $\overline{V}_0 \cap \MP(f)\subset \partial\Delta$, there exists $N \geq 1$, such that $\diam_{\EC}(V_n)<\varepsilon$ for all $n\geq N$, where $V_n$ is any connected component of $f^{-n}(V_0)$.
\end{thm}

\begin{proof}[Proof of Theorem \ref{thm:JLC}]
For simplicity, we denote $f:=f_\lambda$ and $\Delta:=\Delta_\lambda$.
By Lemma \ref{lem:Fatou-finite}, all Fatou components of $f$ are bounded Jordan domains. Hence $J(f)$ is connected. By Lemma \ref{lem:LC-criterion}, it suffices to prove that for any $\varepsilon>0$, there are only finitely many Fatou components of $f$ whose spherical diameters are greater than $\varepsilon$.

For $n\geq 1$, let $U_n$ be a Fatou component such that $f^{\circ n}(U_n)=\Delta$ and $f^{\circ (n-1)}(U_n)\neq\Delta$
Then $f^{\circ n}:U_n\to \Delta$ is conformal and $f^{\circ n}:\partial U_n\to\partial \Delta$ is a homeomorphism.
There exist finitely many bounded Jordan domains $\{V_{0,(i)}:1\leq i\leq M_0\}$ in $\C\setminus\overline{\Delta}$ such that $\overline{V}_{0,(i)}\cap \partial\Delta\neq\emptyset$ is a subarc for all $1\leq i\leq M_0$ and $\partial\Delta\subset\bigcup_{i=1}^{M_0} \overline{V}_{0,(i)}\subset\C\setminus\{0\}$.
Let $\varepsilon>0$ be given. By Theorem \ref{thm:main-pre}, there exists $N_0\geq 1$ such that for all $n\geq N_0$ and $1\leq i\leq M_0$, we have $\diam_{\EC}(V_{n,(i)})<\varepsilon/M_0$, where $V_{n,(i)}$ is any component of $f^{-n}(V_{0,(i)})$. Hence we have $\diam_{\EC}(U_n)<\varepsilon$ for all $n\geq N_0$.

Note that there exists a number $R=R(\varepsilon)>1$ such that if $X\subset\C$ is disjoint with $\D(0,R)$, then $\diam_{\EC}(X)<\varepsilon$.
For any given $n\geq 1$, the Fatou components $U_n$'s of $f$ satisfying $f^{\circ n}(U_n)=\Delta$ and $f^{\circ (n-1)}(U_n)\neq \Delta$ can only accumulate at the infinity.
Hence for given $n\geq 1$, there are only finitely many of $U_n$'s having spherical diameters greater than $\varepsilon$. The proof is complete.
\end{proof}

\section{The round circle for the premodel}\label{sec:model}

In this section, we locate the position of the round circle $\Gamma$ (see Figure \ref{Fig:model}) and prove the {\hyperref[lem:main]{Main Lemma}}.
In fact, we will prove a more result, which implies the {\hyperref[lem:main]{Main Lemma}} immediately.

\subsection{The reduction}

For a Jordan curve $\gamma$ in $\C$, we use $\gamma^{\Ext}$ to denote the component of $\EC\setminus\gamma$ containing $\infty$ and use $\gamma^{\Int}$ to denote the other.
Recall that the critical point and the critical value of $g(z)=ze^z$ are $c:=-1$ and $v:=-1/e$ respectively.

\begin{lem}\label{lem:main-2}
For every $\theta\in(\R/\Z)\setminus\{0\}$, there exist a round circle $\Gamma$ passing through $v=-1/e$ and a Jordan curve $\gamma\subset g^{-1}(\Gamma)$ passing through $c=-1$ such that $0\in\Gamma^{\Ext}\cap\gamma^{\Ext}$, the restriction $g|_\gamma$ is an orientation-preserving homeomorphism onto $\Gamma$, and there is a conformal isomorphism $\phi:\Gamma^{\Ext}\to\gamma^{\Ext}$ fixing $0$ and $\infty$ such that the rotation number satisfies $\rot(g\circ\phi|_\Gamma)=\theta$.
\end{lem}

Here the boundary map $\phi:\Gamma\to\gamma$ is understood as the Carath\'eodory extension of the conformal isomorphism $\phi:\Gamma^{\Ext}\to\gamma^{\Ext}$ and in the following we will still use this convention. Moreover, we will see that $\gamma$ is the boundary of a component of $g^{-1}(\Gamma^{\Int})$.

We first give a decomposition of the dynamical plane of the map $g(z)=ze^z$. For $\ell\in\Z$, we denote
\begin{equation}\label{equ:eta-ell}
\eta_\ell:=\{\eta_\ell(y)=-y\cot y+iy\in\C:\ell\pi<y<(\ell+1)\pi\}.
\end{equation}

\begin{figure}[!htpb]
  \setlength{\unitlength}{1mm}
  \centering
  \includegraphics[width=0.75\textwidth]{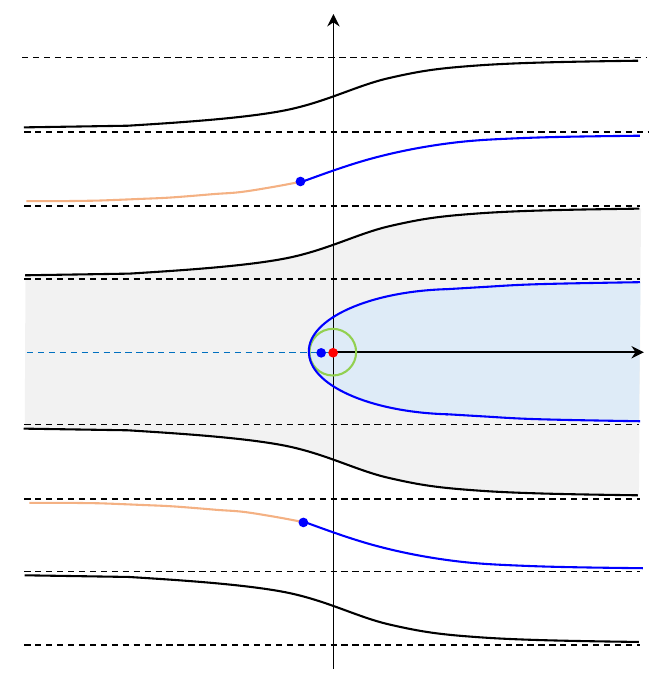}
  \put(-52,3){$-4\pi$}
  \put(-52,15){$-3\pi$}
  \put(-52,26.5){$-2\pi$}
  \put(-52,39){$-\pi$}
  \put(-52,50){$0$}
  \put(-55.3,54.8){\textcolor{blue}{$v$}}
  \put(-62,50){$-1$}
  \put(-49,56){$\D$}
  \put(-52,62){$\pi$}
  \put(-52,73.5){$2\pi$}
  \put(-52,85.5){$3\pi$}
  \put(-52,97.5){$4\pi$}
  \put(-51,106){$y$}
  \put(-3,50){$x$}
  \put(-101,14){$\eta_{-4}$}
  \put(-101,26){$\eta_{-3}$}
  \put(-101,38){$\eta_{-2}$}
  \put(-10,45){$\eta_{-1}$}
  \put(-9,61.5){$\eta_0$}
  \put(-101,68){$\eta_1$}
  \put(-101,80){$\eta_2$}
  \put(-101,92){$\eta_3$}
  \put(-80,57){$V_0$}
  \put(-30,57){\textcolor{blue}{$W_0$}}
  \caption{The preimage of $\R$ under $g(z)=ze^z$ consists of countably many real-analytic curves $\R\cup\{\eta_\ell:\ell\in\Z\}$, where each $\eta_\ell$ is contained in a horizontal strip of hight $\pi$. The black curves are preimages of $[0,+\infty)$ and the blue curves are preimages of $(-\infty,v]$. The domain $V_0$ is bounded by $\eta_1\cup\eta_{-2}\cup[0,+\infty)$ and $W_0$ is bounded by $\eta_0\cup\eta_{-1}\cup\{-1\}$.}
\label{Fig-g-inverse}
\end{figure}

\begin{lem}[see Figure \ref{Fig-g-inverse}]\label{lem:eta-ell}
We have
\begin{equation}\label{equ:R-inverse}
g^{-1}(\R)=\R\cup\bigcup_{\ell\in\Z}\eta_\ell,
\end{equation}
where $\eta_\ell=\overline{\eta_{-\ell-1}}$ is a real-analytic curve and $y\mapsto \re\eta_\ell(y)$ is strictly increasing in $(\ell\pi,(\ell+1)\pi)$ for each $\ell\geq 0$. If $\ell\geq 1$, then $\{z\in\C:\im z=\ell\pi\}$ and $\{z\in\C:\im z=(\ell+1)\pi\}$ are two asymptotes of $\eta_\ell$ and
\begin{equation}\label{equ:asym-1}
\lim_{y\to \ell\pi^+}\re\eta_\ell(y)=-\infty  \text{\quad and\quad}
\lim_{y\to (\ell+1)\pi^-}\re\eta_\ell(y)=+\infty.
\end{equation}
Moreover, $\{z\in\C:\im z=\pi\}$ is the asymptote of $\eta_0$ and
\begin{equation}\label{equ:asym-2}
\lim_{y\to 0^+}\re\eta_0(y)=-1  \text{\quad and\quad}
\lim_{y\to \pi^-}\re\eta_0(y)=+\infty.
\end{equation}
\end{lem}

\begin{proof}
Denote $z=x+iy$, where $x,y\in\R$. Then
\begin{equation}\label{equ:g-z}
g(z)=ze^z=e^x(x\cos y-y\sin y)+i e^x(x\sin y+y\cos y).
\end{equation}
Hence $\im g(z)=0$ is equivalent to $x\sin y+y\cos y=0$.
If $\sin y=0$, then $y=0$. If $\sin y\neq 0$, then
\begin{equation}\label{equ:x-ycoty}
x=-y\cot y, \text{\quad where } y\not\in\pi\Z
\end{equation}
and we have \eqref{equ:R-inverse}.
A direct calculation shows that
\begin{equation}
(-y\cot y)'=\frac{y-\sin y\cos y}{\sin^2 y}>0 \text{\quad for } y>0.
\end{equation}
Hence $\re\eta_\ell(y)$ is strictly increasing in $(\ell\pi,(\ell+1)\pi)$ for each $\ell\geq 0$. By symmetry, we have $\eta_\ell=\overline{\eta_{-\ell-1}}$ and they are real-analytic curves for $\ell\geq 0$.
The statements \eqref{equ:asym-1} and \eqref{equ:asym-2} can be obtained by simple calculations.
\end{proof}

\begin{lem}[see Figure \ref{Fig-g-inverse}]\label{lem:g-inverse}
We have
\begin{equation}\label{equ:g-inverse-1}
g^{-1}([0,+\infty))=[0,+\infty)\cup\bigcup_{k\geq 1}(\eta_{2k-1}\cup \eta_{-2k}).
\end{equation}
Let $V_0$ be the domain bounded by $\eta_1\cup\eta_{-2}\cup[0,+\infty)$. Then
$g:V_0\to\C\setminus[0,+\infty)$
is a proper map of degree two.
\end{lem}

\begin{proof}
By \eqref{equ:g-z}, if $\re g(z)>0$, then equivalently, $z\in \R\cup\bigcup_{\ell\in\Z}\eta_\ell$ and
\begin{equation}\label{equ:x-y-1}
x\cos y-y\sin y>0.
\end{equation}
If $z\in\R$, then $y=0$ and $x>0$. Otherwise,  we have $\sin y\neq 0$ and by \eqref{equ:x-ycoty},
\begin{equation}\label{equ:ysiny}
x\cos y-y\sin y=-y\cot y\cos y-y\sin y=-\frac{y}{\sin y}>0.
\end{equation}
Thus $y\sin y<0$.
If $y>0$, we have $(2k-1)\pi<y<2k\pi$ for some $k\geq 1$, and if $y<0$, we have $-2k\pi<y<(-2k+1)\pi$ for some $k\geq 1$. Thus \eqref{equ:g-inverse-1} holds.
By Proposition \ref{prop:E-preimage}, $g:V_0\to\C\setminus[0,+\infty)$ is a proper map of degree two.
\end{proof}

\begin{lem}[see Figure \ref{Fig-g-inverse}]\label{lem:g-inverse-1}
We have
\begin{equation}\label{equ:g-inverse-2}
g^{-1}((-\infty,v])=\{-1\}\cup\eta_0\cup\eta_{-1}\cup\bigcup_{k\geq 1}(\eta_{2k}'\cup \eta_{-2k-1}'),
\end{equation}
where $\eta_{2k}'=\overline{\eta_{-2k-1}'}=\{\eta_{2k}(y):y\in[y_{2k},(2k+1)\pi)\}$ for some $y_{2k}\in(2k\pi,(2k+1)\pi)$ is a continuous subcurve of $\eta_{2k}$ for each $k\geq 1$.
Let $W_0$ be the domain containing $[0,+\infty)$ which is bounded by $\eta_0\cup\eta_{-1}\cup\{-1\}$. Then $\D\subset W_0\subset\{z\in\C:\re z>-1\}$ and
$g:W_0\to\C\setminus(-\infty,v]$ is conformal.
\end{lem}

\begin{proof}
By \eqref{equ:g-z}, if $\re g(z)\leq v=-1/e$, then equivalently, $z\in \R\cup\bigcup_{\ell\in\Z}\eta_\ell$ and
\begin{equation}\label{equ:x-y-2}
e^x(x\cos y-y\sin y)\leq -1/e.
\end{equation}
If $\sin y=0$, then $y=0$ and $x=-1$. If $\sin y\neq 0$, then $y\not\in \pi\Z$ and $x=-y\cot y$ by \eqref{equ:x-ycoty}. By \eqref{equ:ysiny} and \eqref{equ:x-y-2} we have $e^x\frac{y}{\sin y}\geq \frac{1}{e}$. Hence
\begin{equation}
y\sin y>0 \text{\quad and\quad}
\xi_0(y):=\cos y\frac{y}{\sin y}-\log\frac{y}{\sin y}-1\leq 0.
\end{equation}
If $y>0$, we have $2k\pi<y<(2k+1)\pi$ for some $k\geq 0$, and if $y<0$, we have $(-2k-1)\pi<y<-2k\pi$ for some $k\geq 0$.

By a direct calculation, we have
\begin{equation}
\xi_0'(y)=-\frac{(y-\sin y\cos y)^2+\sin^4 y}{y\sin^2 y}<0 \text{\quad for } y>0.
\end{equation}
Hence $\xi_0(y)$ is strictly decreasing in $(2k\pi,(2k+1)\pi)$ for each $k\geq 0$. Note that
\begin{equation}\label{equ:asym-3}
\lim_{y\to 2k\pi^+}\xi_0(y)=+\infty  \text{\quad and\quad}
\lim_{y\to (2k+1)\pi^-}\xi_0(y)=-\infty
\end{equation}
for each $k\geq 1$ and
\begin{equation}\label{equ:asym-4}
\lim_{y\to 0^+}\xi_0(y)=0  \text{\quad and\quad}
\lim_{y\to \pi^-}\xi_0(y)=-\infty.
\end{equation}
Hence there exists a unique $y_{2k}\in(2k\pi,(2k+1)\pi)$ such that $\xi_0(y)\leq 0$ for all $y\in[y_{2k},(2k+1)\pi)$, where $k\geq 1$. Moreover, $\xi_0(y)\leq 0$ for all $y\in(0,\pi)$.
By symmetry, we have $\eta_{2k}'=\overline{\eta_{-2k-1}'}$ and it is a continuous subcurve of $\eta_{2k}$ for all $k\geq 1$. Thus \eqref{equ:g-inverse-2} holds.

Let $W_0$ be the domain containing $[0,+\infty)$ which is bounded by $\eta_0\cup\eta_{-1}\cup\{-1\}$ (see Figure \ref{Fig-g-inverse}). Then
\begin{equation}\label{equ:W-0}
W_0=\{x+iy:|y|<\pi \text{ and }x>-y\cot y\}.
\end{equation}
In Lemma \ref{lem:eta-ell}, we have proved that $y\mapsto -y\cot y$ is strictly increasing in $(\ell\pi,(\ell+1)\pi)$ for each $\ell\geq 0$. Hence $W_0\subset\{z\in\C:\re z>-1\}$.
To prove $\D\subset W_0$, it suffices to verify that $-y\cot y\leq-\sqrt{1-y^2}$ for all $0<y\leq 1$.
This reduces to the inequality $y\geq \sin y$, which holds for $0<y\leq 1$.
Finally, by Proposition \ref{prop:E-preimage}, $g:W_0\to\C\setminus(-\infty,v]$ is conformal.
\end{proof}

By Lemma \ref{lem:g-inverse-1}, the map $g$ has an inverse branch
\begin{equation}\label{equ:h-0}
h_0:\C\setminus (-\infty,v]\to W_0.
\end{equation}
Moreover, if $D$ is Jordan domain in $\C\setminus[0,+\infty)$ with $v\in\partial D$ and $0\not\in\partial D$, then by Lemma \ref{lem:g-inverse}, $g^{-1}(D)\cap V_0$ consists of two Jordan domains which touch only at $c=-1$.

For a closed Jordan domain $K\subset\C$, there exists a unique conformal map
\begin{equation}\label{equ:psi-K}
\psi_K:\EC\setminus K\to\EC\setminus\overline{\D}
\end{equation}
such that $\psi_K(\infty)=\infty$ and $\psi_K'(\infty)>0$. It is well-known that the Green function
\begin{equation}\label{equ:G-K}
G_K(z):=\log|\psi_K(z)|
\end{equation}
is positive and harmonic in $\C\setminus K$, and moreover, it can be extended continuously to the boundary and equal zero there. Furthermore, if $K_1\subset K_2$ are two closed Jordan domains in $\C$, then $G_{K_1}\ge G_{K_2}$ outside $K_2$.

\begin{figure}[!htpb]
  \setlength{\unitlength}{1mm}
  \centering
  \includegraphics[width=0.75\textwidth]{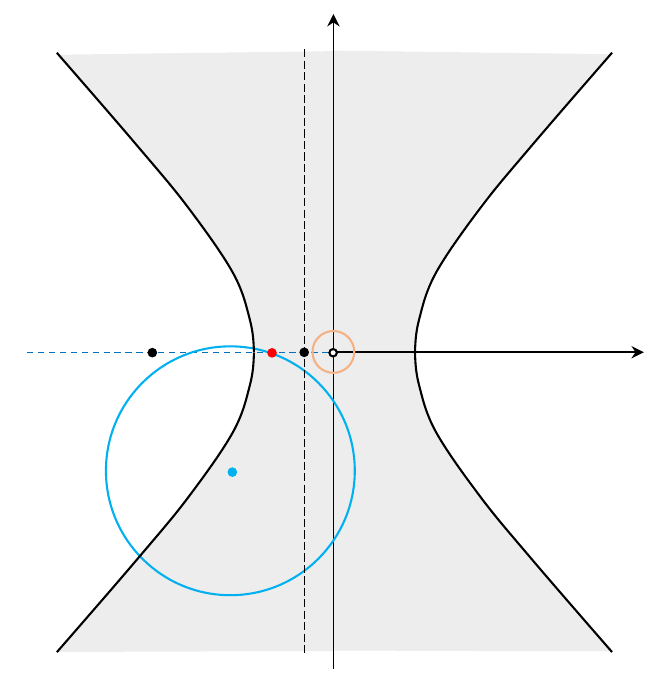}
  \put(-73.5,29.5){$v+t^2$}
  \put(-68,10.5){$\Gamma_t$}
  \put(-52,50){$0$}
  \put(-60.3,56){$\frac{v}{2}$}
  \put(-64,50){$v$}
  \put(-49,56){$\T_r$}
  \put(-43,83){$\Omega$}
  \put(-28,75){$x^2-y^2=-\frac{v}{2}$}
  \put(-51,106){$y$}
  \put(-3,50){$x$}
  \put(-83,56){$c=-1$}
  \put(-39,49){$\sqrt{-\frac{v}{2}}$}
  \caption{The parameter domain $\Omega$ (colored gray) for the family of circles $\{\Gamma_t\}_{t\in\Omega\setminus\{0\}}$, which is bounded by two hyperbolas. The center of $\Gamma_t$ lies in the left half plane $\{t\in\C:\re t<\frac{v}{2}\}$. A special circle $\Gamma_t$ is drawn. It passes through the critical value $v=-1/e$ and is disjoint with the non-negative real axis.}
\label{Fig-Gamma-Omega}
\end{figure}

\medskip
To prove Lemma \ref{lem:main-2}, we shall construct a family of round circles $\{\Gamma_t\}_{t\in\Omega\setminus\{0\}}$ such that $v\in\Gamma_t$ and $\overline{\Gamma_t^{\Int}}\subset \C\setminus[0,+\infty)$ for all $t\in\Omega\setminus\{0\}$. We prove the existence of $t_0\in\Omega\setminus\{0\}$ such that Lemma \ref{lem:main-2} holds for the circle $\Gamma=\Gamma_{t_0}$.
Specifically, we set
\begin{equation}\label{equ:Gamma-t}
\Gamma_t:=\partial\D(v+t^2,|t|^2),
\end{equation}
where
\begin{equation}\label{equ:Omega}
t\in\Omega\setminus\{0\} \text{ and } \Omega:=\{t\in\C:\re(t^2)<-\tfrac{v}{2}\}.
\end{equation}
We will see later that it is convenient to choose $t^2$ but not $t$ in the parametrization.
Note that $v\in\Gamma_t$. Since $\re(v+t^2)<\frac{v}{2}$, we have
\begin{equation}\label{equ:Gamma-disjoint}
\Gamma_t\cap [0,+\infty)=\emptyset \text{\quad for all } t\in\Omega\setminus\{0\}.
\end{equation}
See Figure \ref{Fig-Gamma-Omega}.

\begin{figure}[!htpb]
  \setlength{\unitlength}{1mm}
  \centering
  \includegraphics[width=0.95\textwidth]{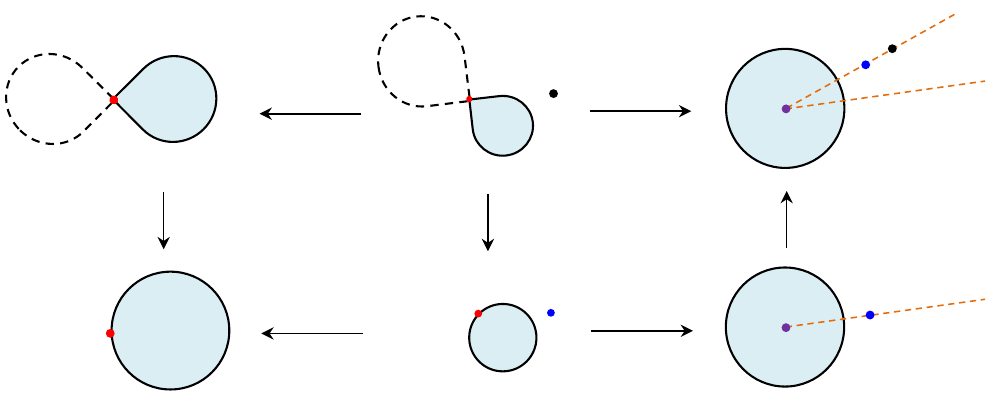}
  \put(-110,24.5){$q(z)=w^2$}
  \put(-67,24.5){$g(z)=ze^z$}
  \put(-74,2){$\Gamma_t$}
  \put(-61,8){$0$}
  \put(-73,12){$v$}
  \put(-123,8){$0$}
  \put(-104,2){$\D(1,1)$}
  \put(-120.3,36.5){$0$}
  \put(-114,40){$K_+$}
  \put(-130,40){$K_-$}
  \put(-95,12){$\frac{z-v}{t^2}$}
  \put(-97,43){$\frac{\chi(z+1)}{\sqrt{2e}\,t}$}
  \put(-60.8,37.8){$0$}
  \put(-75,34){$\gamma_t$}
  \put(-71,49){$\gamma_{-t}$}
  \put(-50,42){$\psi_t^\gamma$}
  \put(-50,12){$\psi_t^\Gamma$}
  \put(-74.5,42){$c$}
  \put(-29,36){$0$}
  \put(-29,6){$0$}
  \put(-20,3){$\D$}
  \put(-20,32){$\D$}
  \put(-27,23){$R_{s_t}$}
  \put(-20,14.5){\small{$\psi_t^{\Gamma}(0)$}}
  \put(-19,50.5){\small{$\psi_t^{\gamma}(0)$}}
  \put(-19,41.7){\rotatebox{28}{\small{$R_{s_t}\circ\psi_t^{\Gamma}(0)$}}}
  \caption{Some domains and maps which are used in the proof of Lemma \ref{lem:main-2}.}
\label{Fig-G-rho}
\end{figure}

By Lemma \ref{lem:g-inverse} and \eqref{equ:Gamma-disjoint}, for $t\in \Omega\setminus\{0\}$, there exists exactly two bounded Jordan curves $\gamma_t$ and $\gamma_{-t}$ in $V_0$ with $c\in\gamma_t\cap\gamma_{-t}$ such that $\gamma_t$ and $\gamma_{-t}$ move continuously about $t$, and
\begin{equation}\label{equ:gamma-t-pm}
g:\gamma_t^{\Int}\to \Gamma_t^{\Int} \text{\quad and\quad} g:\gamma_{-t}^{\Int}\to \Gamma_t^{\Int}
\end{equation}
are conformal. In fact, $\gamma_t$ and $\gamma_{-t}$ play equal roles in the rest arguments.
Let
\begin{equation}
\psi_t^\Gamma:=\psi_{\overline{\Gamma_t^{\Int}}}:\Gamma_t^{\Ext}\to\EC\setminus\overline{\D}
 \text{\quad and\quad}
 \psi_t^\gamma:=\psi_{\overline{\gamma_t^{\Int}}}:\gamma_t^{\Ext}\to\EC\setminus\overline{\D}
\end{equation}
be the normalized conformal maps as in \eqref{equ:psi-K}.
In particular,
\begin{equation}\label{equ:psi-t-Gamma}
\psi_t^\Gamma(z)=(z-v-t^2)/|t|^2.
\end{equation}
Denote the rigid rotation
\begin{equation}\label{equ:R-s-t}
R_{s_t}(\zeta):=e^{2\pi i s_t}\zeta=\frac{\psi_t^\gamma(0)}{|\psi_t^\gamma(0)|}/\frac{\psi_t^\Gamma(0)}{|\psi_t^\Gamma(0)|}\zeta.
\end{equation}
Then
\begin{equation}\label{equ:phi-t}
\phi_t:=(\psi_t^\gamma)^{-1}\circ R_{s_t}\circ\psi_t^\Gamma:\Gamma_t^{\Ext}\to\gamma_t^{\Ext}
\end{equation}
is a conformal map fixing the infinity. Moreover, $R_{s_t}\circ\psi_t^\Gamma(0)$ and $\psi_t^\gamma(0)$ have the same argument. See Figure \ref{Fig-G-rho}.

For $t\in\Omega\setminus\{0\}$, define
\begin{equation}\label{equ:G-t}
 G(t):=\log|\psi_t^\Gamma(0)|-\log|\psi_t^\gamma(0)|.
\end{equation}
Then we have
\begin{equation}\label{equ:phi-t-G-equi}
\phi_t(0)=0\iff G(t)=0.
\end{equation}
For the orientation-preserving homeomorphism
\begin{equation}\label{equ:g-t}
g_t:=\psi_t^\Gamma\circ g\circ(\psi_t^\gamma)^{-1}\circ R_{s_t}:\T\to\T,
\end{equation}
we denote the rotation number of $g_t|_{\T}$ by
\begin{equation}
\rho(t):=\rot(g_t)\in\R/\Z.
\end{equation}
The quantities $G(t)$ and $\rho(t)$ are continuous with respect to $t\in\Omega\setminus\{0\}$ (the continuity of $\psi_t^\Gamma$ and $\psi_t^{\gamma}$ follows by Carath\'{e}odory Kernal Theorem, and the continuity of the rotation number follows by \cite[Proposition 11.1.6, p.\,390]{KH95}).
Note that $g_t:\T\to\T$ is conjugate to $g\circ\phi_t:\Gamma_t\to\Gamma_t$.
To prove Lemma \ref{lem:main-2}, it suffices to find a parameter $t_0\in\Omega\setminus\{0\}$ such that $G(t_0)=0$ and $\rho(t_0)=\theta$.

\subsection{The sufficiently small and large parameters}

In this subsection, we consider the values of $G(t)$ and $\rho(t)$ when $t\in\Omega\setminus\{0\}$ is small (resp., large) enough.

\begin{prop}\label{prop:small-circle}
There exists $0<r_0<|v|/4$ such that $G(t)>0$ for all $0<|t|\leq r_0$ and the winding number of $e^{2\pi i\rho(t)}$ around $0$ is $1$ as the parameter $t$ makes one full turn around the circle $\T_{r}$ for any $r\in(0,r_0]$.
\end{prop}

\begin{proof}
We first construct a local model of $g(z)=ze^z$ near the critical point $c=-1$ by considering the square mapping $q(w)=w^2$. The map $q$ maps two closed Jordan domains
\begin{equation}
K_\pm:=\{w\in\C:|w^2-1|\le1,\ \pm\re w\ge0\}
\end{equation}
homeomorphically onto the closed disk $\overline{\D}(1,1)$.
Note that $g$ is local $2$-to-$1$ at $c$ and has the Taylor expansion
\begin{equation}
g(z)=v+\tfrac{1}{2e}(z-c)^2+O((z-c)^3).
\end{equation}
There exists a conformal map $\widetilde{\chi}_t(z):=\frac{\chi(z+1)}{\sqrt{2e}\,t}$ defined in a neighborhood of $\overline{\gamma_t^{\Int}}\cup\overline{\gamma_{-t}^{\Int}}$ to that of $K_+\cup K_-$ which is the lift of $\frac{z-v}{t^2}$ under $g$ and $q$ such that
\begin{equation}\label{equ:chi}
\chi(0)=0,~~\chi'(0)=1 \text{ and } \widetilde{\chi}_t(\overline{\gamma_{\pm t}^{\Int}})=K_\pm.
\end{equation}
Note that $\chi$ is defined in a fixed neighborhood of $0$, which is independent of small $t$.
See Figure \ref{Fig-G-rho}.

Denote $K_t:=\{\frac{z+1}{\sqrt{2e}\,t}:z\in \overline{\gamma_t^{\Int}}\}$. Let
\begin{equation}
\psi_{K_t}:\EC\setminus K_t\to\EC\setminus\overline{\D} \text{\quad and\quad} \psi_{K_+}:\EC\setminus K_+\to\EC\setminus\overline{\D}
\end{equation}
be the normalized conformal mappings defined as in \eqref{equ:psi-K}. We have
\begin{equation}\label{equ:psi-t-c}
\psi_t^\gamma(z)=\tfrac{t}{|t|}\,\psi_{K_t}\big(\tfrac{z+1}{\sqrt{2e}\,t}\big) \text{\quad and\quad} \psi_t^\gamma(c)=\tfrac{t}{|t|}\,\psi_{K_t}(0)\in\T.
\end{equation}
By \eqref{equ:chi}, the compact set $K_t$ converges to $K_+$ in the kernal sense as $t\to 0$.
Moreover, there exists a constant $M>1$ which is independent of small $t$ such that $\partial K_t$ is an $M$-quasi-circle for all small $t\neq 0$. Then according to \cite[Theorem 9.11, p.\,283]{Pom75}, $\psi_{K_t}^{-1}$ converges to $\psi_{K_+}^{-1}$ uniformly on $\EC\setminus\D$ as $t\to 0$. In particular, we have
\begin{equation}\label{equ:psi-limit}
\psi_{K_t}(0)\to \psi_{K_+}(0) \text{ as } t\to 0.
\end{equation}
A direct calculation shows that
\begin{equation}\label{equ:psi-t-v}
\psi_t^\Gamma(z)=\tfrac{z-v-t^2}{|t|^2} \text{\quad and\quad} \psi_t^\Gamma(v)=-\tfrac{t^2}{|t|^2}\in\T.
\end{equation}
Hence by \eqref{equ:psi-t-c}, \eqref{equ:psi-limit} and \eqref{equ:psi-t-v}, there exists a small $r_1>0$ such that as the parameter $t=re^{i\alpha}$ makes one full turn around the circle $\T_{r}$ for any $r\in(0,r_1]$, then the winding numbers of $\psi_t^\gamma(c)=e^{i\alpha}\,\psi_{K_t}(0)$ and $\psi_t^\Gamma(v)=-e^{2i\alpha}$ around $0$ is $1$ and $2$ respectively.

By \eqref{equ:psi-t-c} and \eqref{equ:psi-t-v}, we have 
\begin{equation}\label{equ:psi-t-2}
\psi_t^\gamma(0)=\tfrac{t}{|t|}\,\psi_{K_t}\big(\tfrac{1}{\sqrt{2e}\,t}\big) \text{\quad and\quad}
\psi_t^\Gamma(0)=\tfrac{1}{e|t|^2}-\tfrac{t^2}{|t|^2}.
\end{equation} In $t\neq 0$ is small enough, then $\arg\psi_t^\gamma(0)$ and $\arg\psi_t^\Gamma(0)$ are very small.
Hence there exists a small $r_2>0$ such that as the parameter $t$ makes one full turn around the circle $\T_{r}$ for any $r\in(0,r_2]$, then the winding number of
$e^{2\pi i s_t}=\frac{\psi_t^\gamma(0)}{|\psi_t^\gamma(0)|}/\frac{\psi_t^\Gamma(0)}{|\psi_t^\Gamma(0)|}$
around $0$ is $0$. By \eqref{equ:g-t}, we have
\begin{equation}
g_t\big(R_{s_t}^{-1}\circ \psi_t^\gamma(c)\big)=\psi_t^\Gamma(v)\in\T.
\end{equation}
This implies that the winding number of $e^{2\pi i\rho(t)}$ around $0$ is $1$ as the parameter $t=re^{i\alpha}$ makes one full turn around $\T_{r}$ for any $r\in(0,r_0']$ with $r_0':=\min\{r_1,r_2\}$.

\medskip
To verify the sign of $G(t)$, note that $\gamma_t^{\Int}$ contains a disk of radius $a_0|t|$ whose center tends to $c$, and is contained in $\D(c,A_0|t|)$, where $a_0,A_0>0$ are independent of small $t$.
This implies that (see \eqref{equ:psi-t-2})
\[
 \log|\psi_t^\gamma(0)|=\log\tfrac{1}{|t|}+O(1) \text{\quad and\quad}
 \log|\psi_t^\Gamma(0)|=2\log\tfrac{1}{|t|}+O(1).
\]
Thus there exists a small $r_3>0$ such that $G(t)=\log\tfrac{1}{|t|}+O(1)>0$ for all $0<|t|\leq r_3$. The proof is complete by setting $r_0:=\min\{r_0',r_3, |v|/4\}>0$.
\end{proof}

To estimate the value of $G(t)$ for large $t$, we will apply the following \textit{Harnack's principle}. 

\begin{lem}[Harnack's principle]\label{lem:harnack}
If $u>0$ is harmonic in $\D(a,R)$, then for any $|z-a|=r<R$,
\begin{equation}\label{equ:harmonic-disk}
\frac{R-r}{R+r}u(a)\le u(z)\le\frac{R+r}{R-r}u(a).
\end{equation}
Moreover, if $u>0$ is harmonic in $\{z\in\C:\re z>0\}$, then for $0<x_1\le x_2$,
\begin{equation}\label{equ:harmonic-plane}
\frac {x_1}{x_2} u(x_2)\le u(x_1)\le\frac{x_2}{x_1} u(x_2).
\end{equation}
\end{lem}

The proof of the first statement of Lemma \ref{lem:harnack} is based on using the Poisson formula (see \cite[p.\,243]{Ahl78}) and the second statement can be obtained by transforming the right half plane to the unit disk under $z\mapsto(z-x_2)/(z+x_2)$.
Note that one can apply \eqref{equ:harmonic-disk} successively along a fixed finite chain of disks in a domain to compare the values at the two ends, and the corresponding comparisons to \eqref{equ:harmonic-plane} hold also in translated or rotated half-planes.

\medskip
For a closed Jordan domain $K\subset\C$, recall that $G_K$ is the Green function in $\C\setminus K$ defined in \eqref{equ:G-K}.

\begin{lem}\label{lem:est-green}
Let $K\subset\overline{\D}_L$ be a closed Jordan domain, where $L\geq 1$.
\begin{enumerate}
\item
if $K\subset\{z\in\C: \im z\ge 0\}\setminus\{0\}$, then
\begin{equation}
G_K(0)\ge\frac{d\log2}{12L}, \text{\quad where }d=\min\{1,\dist_{\C}(0,K)\};
\end{equation}

\item if $K\subset\{z\in\C:\im z\ge-\eta\}$, where $0\le\eta<s/2$ with $0<s\le1$, then
\begin{equation}
G_K(-is)\ge \frac{s\log2}{4L};
\end{equation}

\item if $u>0$ is harmonic in $\{z\in\C: |z|<3s \text{ and } \re z>0\}$, then
\begin{equation}\label{eq:est-halfdisc}
u(x)\ge\frac{x}{2s-x}u(s)\ge\frac{x}{2s}u(s), \text{\quad where } 0<x<\frac{s}{2}.
\end{equation}
\end{enumerate}
\end{lem}

\begin{proof}
(a) Since $K\subset\overline{\D}_L$, we have $G_K(-2iL)\ge\log2$. By the half-plane Harnack inequality \eqref{equ:harmonic-plane}, we have
\[
G_K(-iy)\ge\frac{y}{2L}\log2, \text{\quad for all } 0<y\le2L.
\]
Take $y=\frac{d}{2}$ and compare $G_K(0)$ with $G_K(-i\frac{d}{2})$ in $\D_d$, by the disk Harnack inequality \eqref{equ:harmonic-disk} we obtain
$G_K(0)\ge G_K(-i\frac{d}{2})/3\ge d\log2/(12L)$.

(b) Since $K\subset\overline{\D}_L$, we have $G_K(-i(2L+\eta))\ge\log2$.
Compare this point with $-is$ in $\{z\in\C:\im z<-\eta\}$. Their distances from the boundary are $2L$ and $s-\eta$ respectively. By the half-plane Harnack inequality \eqref{equ:harmonic-plane} we have
\[
G_K(-is)\ge\frac{s-\eta}{2L}\log2\ge\frac{s\log2}{4L}.
\]

(c) Note that the disk $\D(s,s)$ lies in the given right half-disk, and the distance from $x$ to the center of this disk is $s-x$. Thus \eqref{eq:est-halfdisc} holds by \eqref{equ:harmonic-disk}.
\end{proof}


Let $A$ and $B$ be two non-negative numbers. In the following, we use $A\preceq B$ (resp., $A\succeq B$) to denote $A\le cB$ (resp., $B\le cA$) for a constant $c>0$ which is independent of the varying large parameter $t$, and $A\asymp B$ to mean comparable, i.e., $A\preceq B$ and $B\preceq A$. Denote $\log^+x:=\max\{\log x,0\}$.

\begin{prop}\label{prop:est-infinity}
There exists $R_0>1$ such that $G(t)<0$ for all $t\in\Omega$ with $|t|\geq R_0$.
\end{prop}

\begin{proof}
Let $C=a+ib$ be the center of the disk
\begin{equation}
D:=\Gamma_t^{\Int}=\D(v+t^2,|t|^2)\subset\C\setminus[0,+\infty),
\end{equation}
where $a=v+\re(t^2)=-\frac{1}{e}+\re(t^2)$ and $b=\im(t^2)$ (see \eqref{equ:Gamma-disjoint}). By \eqref{equ:Omega}, $\re(t^2)<-\frac{v}{2}$. Then
\begin{equation}\label{equ:a0c0}
a<\tfrac{v}{2} \text{\quad and\quad} \kappa:=v(2a-v)>0.
\end{equation}
Denote
\begin{equation}\label{equ:R-r-delta}
R:=|C|=|v+t^2|,\quad r:=|t|^2 \text{\quad and\quad} \delta=R-r=\frac{\kappa}{R+r}>0.
\end{equation}
Thus by \eqref{equ:psi-t-Gamma} we have
\begin{equation}\label{equ:Theta-1}
G_{\overline{D}}(0)=\log\frac{R}{r}=\log\Big(1+\frac{\delta}{r}\Big)=\frac{1}{2}\log\Big(1+\frac{\kappa}{r^2}\Big).
\end{equation}
Let $K$ be either $\overline{\gamma_t^{\Int}}$ or $\overline{\gamma_{-t}^{\Int}}$ in $V_0$ (see \eqref{equ:gamma-t-pm}), where $V_0$ is the domain introduced in Lemma \ref{lem:g-inverse}. Then $|\im z|<2\pi$ for all $z\in K$. Hence if $\delta\le|g(z)|=|ze^z|\le M$ for $z\in K$, then
\begin{equation}\label{eq:est-rectangle}
 |z|\preceq 1+\log^+ M+\log^+(1/\delta).
\end{equation}
It suffices to prove that $G_{\overline{D}}(0)/G_K(0)<1$ for large $t\in\Omega$ (see \eqref{equ:G-t}). We divide the arguments into two cases.

\medskip

\emph{Case 1: $a\le v$.}
Note that $D\cap\R=(2a-v,v)$ (which is empty when $a=v$) and $g(\mathbb R)=[v,\infty)$.
Thus the interior of $K$ is disjoint with $\R$ and must be contained in upper or lower half-plane.
Since $\delta=R-r\le |v|=\frac{1}{e}$, we have $R\asymp r=|t|^2$ for large $t$. By \eqref{equ:a0c0} and \eqref{equ:R-r-delta}, we have
\begin{equation}
\kappa\ge v^2=1/e^2 \text{\quad and\quad} R^{-1}\preceq\delta\le |v|=1/e.
\end{equation}
Since $|g(z)|\le e|z|$ for $|z|\le1$ and from $\delta\le |v|<1$, by \eqref{equ:R-r-delta} and \eqref{eq:est-rectangle} we have
\begin{equation}
 K\subset \overline{\D}_L \text{\quad with }L\preceq \log R \text{\quad and\quad} \dist_{\C}(0,K)\ge\delta/e.
\end{equation}
By \eqref{equ:Theta-1} and Lemma \ref{lem:est-green}(a), we have $G_{\overline{D}}(0)\preceq\delta/R$ and
\begin{equation}\label{eq:est-firstcase}
 \frac{G_{\overline{D}}(0)}{G_K(0)}\preceq\frac{\log R}{R}\to 0 \quad \text{as }t\to\infty.
\end{equation}
Thus $G_{\overline{D}}(0)/G_K(0)<1$ for large $t\in\Omega$.

\medskip
\emph{Case 2: $v<a<v/2$.}
Note that $D\cap\R=(v,2a-v)\subset(v,0)$ and $0<\kappa\le v^2$. By \eqref{equ:Theta-1} we have
\begin{equation}\label{eq:est-target-small}
G_{\overline{D}}(0)\preceq\kappa/R^2.
\end{equation}
Making a complex conjugation if necessary, we assume that $b>0$, where $b=\im(t^2)$ is the imaginary of the center of the circle $\Gamma_t$.
Then
\begin{equation}\label{equ:imag-1}
D\subset\{z\in\C:\im z>-\big(|t|^2-\im(t^2)\big)\}.
\end{equation}
Since $0<\re(t^2)<|v|/2$, we have
\begin{equation}\label{equ:imag-2}
0<|t|^2-\im(t^2)\to 0 \text{\quad as }t\to\infty.
\end{equation}
By Lemma \ref{lem:g-inverse-1} and \eqref{equ:h-0}, $K$ is either contained in $W_0\cup\{-1\}=h_0(\C\setminus (-\infty,v])\cup\{-1\}$ or contained in $(V_0\setminus\overline{W}_0)\cup\{-1\}$.
Thus there exist $A_0>0$ and $0<\eta_R\to0$ as $R\to\infty$ such that by \eqref{eq:est-rectangle}, \eqref{equ:imag-1} and \eqref{equ:imag-2}, we have
\begin{itemize}
\item[(i)] if $K\subset W_0\cup\{-1\}$, then
\begin{equation}\label{eq:est-principal-position}
 K\subset\{z\in\C:\im z\ge-\eta_R\}\cap \D(0,L_1), \text{ where } L_1=A_0\log R;
\end{equation}
\item[(ii)] if $K\subset (V_0\setminus\overline{W}_0)\cup\{-1\}$, then
\begin{equation}\label{eq:est-nonprincipal-position}
 K\subset\{z\in\C:\im z\le \eta_R\}\cap \D(0,L_2),
\text{ where } L_2=A_0\bigl(\log R+\log\tfrac{1}{\kappa}\bigr).
\end{equation}
\end{itemize}

\textit{We first consider (ii)}. In this case $K\cap\D=\emptyset$ by Lemma \ref{lem:g-inverse-1}.
Fixing $0<s<1/4$, the disk Harnack inequality \eqref{equ:harmonic-disk} and Lemma \ref{lem:est-green}(b) therefore give
\begin{equation}
G_K(0)\succeq G_K(is)\succeq \frac1{\log R+\log(1/\kappa)}.
\end{equation}
Since $0<\kappa\le v^2$, we have
\begin{equation}
 \frac{G_{\overline{D}}(0)}{G_K(0)}\preceq \frac{\kappa(\log R+\log(1/\kappa))}{R^2}\preceq \frac{\log R+1}{R^2}\to 0.
\end{equation}
Thus $G_{\overline{D}}(0)/G_K(0)<1$ for large $t\in\Omega$.

\medskip
\textit{Finally we consider (i)}.
In a fixed small neighborhood $\mathcal{N}=\D_{\varepsilon_0}$ of $0$, we define
\begin{equation}\label{eq:est-Q}
\widetilde{g}_C(z):=R\,\Log\left(1-\frac{g(z)}{C}\right)+R\log\frac Rr.
\end{equation}
Here $\Log$ denotes the analytic branch near $1$ normalized by $\Log 1=0$.
Since $C=a+bi=v+t^2$ with $v<a=\re C<v/2$ and $b>0$, by \eqref{equ:R-r-delta} we have
\begin{equation}\label{equ:R-log}
-\frac{R}{C}=-\frac{|C|}{C}\to i \text{\quad and\quad} R\log\frac{R}{r}=\frac{R}{2}\log\Big(1+\frac{\kappa}{r^2}\Big)\preceq\frac{1}{R}.
\end{equation}
Thus $\widetilde{g}_C\to ig$ locally uniformly as $R\to\infty$. Shrinking $\mathcal{N}$ if necessary, we assume that $g$ and $\widetilde{g}_C$ are univalent in $\mathcal{N}$ and $\widetilde{g}_C(\mathcal{N})$ contains a fixed disk $\D_{4s}$ for large $C$, where $s>0$ is small.
The inverse $\widetilde{h}_C:\D_{4s}\to\mathcal{N}$ of $\widetilde{g}_C$ satisfies $\widetilde{h}_C(\zeta)\to h_0(-i\zeta)$ as $C\to\infty$, where $h_0$ is an inverse branch of $g$ defined in \eqref{equ:h-0}.
By \eqref{eq:est-Q} and \eqref{equ:R-log}, we have
\begin{equation}
 \re \widetilde{g}_C(z)=R\log\frac{|g(z)-C|}{r} \text{\quad and\quad} \widetilde{g}_C(0)=R\log\frac{R}{r}\preceq\frac{1}{R}.
\end{equation}

Note that $\{z\in\mathcal{N}:\re \widetilde{g}_C(z)>0\}\cap K=\{z\in\mathcal{N}:|g(z)-C|>r\}\cap K=\emptyset$.
Consequently, $u:=G_K\circ \widetilde{h}_C$ is positive and harmonic in the fixed right half-disk $\{\zeta\in\C:\re\zeta>0 \text{ and }|\zeta|<3s\}$.
Set $x_C:=\widetilde{g}_C(0)=RG_{\overline{D}}(0)\preceq R^{-1}$. Increasing $R$ if necessary, we assume that $0<x_C<s/2$. By Lemma \ref{lem:est-green}(c), we have
\[
G_K(0)=u(x_C)\geq\frac{x_C}{2s}u(s)=\frac{RG_{\overline{D}}(0)}{2s}G_K(\widetilde{h}_C(s)).
\]
Since $\widetilde{h}_C(s)\to h_0(-is)$, the point $\widetilde{h}_C(s)$ stays in a fixed compact subset of the lower half-plane for large $R$.
By \eqref{eq:est-principal-position} and Lemma \ref{lem:est-green}(b), we have $G_K(\widetilde{h}_C(s))\succeq1/\log R$.
This implies that $G_K(0)/G_{\overline{D}}(0)\succeq R/\log R$ and hence $G_{\overline{D}}(0)/G_K(0)<1$ for large $R$.
The proof is complete.
\end{proof}

\subsection{The parameter on the finite boundary}

In this subsection, we consider the values of $G(t)$ and $\rho(t)$ as $t$ tends to the finite point on the boundary $\partial\Omega$.

\begin{prop}\label{prop:est-finite-boundary}
If $t_n\in \Omega\setminus\{0\}$ tends to a finite boundary point $t_*\in\partial\Omega$ and $G(t_n)\ge0$ for all $n$,
then $\rho(t_n)\to0$ in $\mathbb R/\mathbb Z$.
\end{prop}

\begin{proof}
Let $C_n=a_n+ib_n$ be the center of the disk $\Gamma_{t_n}^{\Int}=\D(v+t_n^2,|t_n|^2)\subset\C\setminus[0,+\infty)$, where $a_n=v+\re(t_n^2)$ and $b=\im(t_n^2)$. By \eqref{equ:Omega}, $\re(t_n^2)<-\frac{v}{2}$. Then
\begin{equation}\label{equ:ancn}
a_n<\tfrac{v}{2} \text{\quad and\quad} \kappa_n:=v(2a_n-v)>0.
\end{equation}
Denote $r_n:=|t_n|^2$, $D_n:=\Gamma_{t_n}^{\Int}$ and $K_n:=K_{t_n}$, where $K_{t_n}$ is either $\overline{\gamma_{t_n}^{\Int}}$ or $\overline{\gamma_{-t_n}^{\Int}}$ in $V_0$ (see \eqref{equ:gamma-t-pm}). If $t_n\to t_*\in\partial\Omega$, then
\begin{equation}
C_n\to C_*,\quad \re C_*=\tfrac{v}{2},\quad r_n\to r_*=|C_*|\ge \tfrac{|v|}{2} \text{\quad and\quad}\kappa_n\to0^+.
\end{equation}
By \eqref{equ:Theta-1}, we have
\begin{equation}\label{eq:est-finite-target}
G_{\overline{D}_n}(0)\preceq \kappa_n.
\end{equation}

Since $a_n=\re C_n\to \frac{v}{2}>v$, without loss of generality, we assume that $v<a_n<\frac{v}{2}$ for all $n$. Then $D_n\cap\R=(v,2a_n-v)\subset(v,0)$.
Note that the disks $\{D_n\}$ are uniformly bounded (the bound depends on $t_*$).
By \eqref{eq:est-rectangle}, there are fixed $A\ge1$ and $M>1$ (depending on $t_*$) such that
\[
 K_n\subset[-L_n,M]+i[-2\pi,2\pi], \text{\quad where }
 L_n=A\bigl(1+\log\tfrac{1}{\kappa_n}\bigr).
\]
The above rectangle is contained in the disk $\D(-T_n,T_n+M+\tfrac12)$, where $T_n=\max\{L_n,4\pi^2+1\}$, since $|z+T_n|^2\le(T_n+M)^2+4\pi^2<(T_n+M+\tfrac12)^2$ for each point $z$ in this rectangle.
At the point $x_0:=M+1$, domain monotonicity gives
\begin{equation}\label{eq:est-finite-comparison}
 G_{K_n}(x_0)\ge\log\frac{T_n+M+1}{T_n+M+1/2}\succeq \frac{1}{L_n}.
\end{equation}

By Lemma \ref{lem:g-inverse-1} and \eqref{equ:h-0}, $K_n$ is either contained in $W_0\cup\{-1\}=h_0(\C\setminus (-\infty,v])\cup\{-1\}$ or contained in $(V_0\setminus\overline{W}_0)\cup\{-1\}$.
If there are infinitely many $n$'s such that $K_n\subset (V_0\setminus\overline{W}_0)\cup\{-1\}$, without causing confusion, we still use the subscript ``$n$" to represent this subsequence.
Note that
\[
\mathcal H(z)=\sqrt{1+eg(z)}:W_0\to\{\zeta\in\C:\re\zeta>0\}
\]
is a conformal isomorphism, $\mathcal H(0)=1$ and $\mathcal H(x_0)>1$ is a fixed constant.
Thus $G_{K_n}\circ\mathcal H^{-1}$ is positive and harmonic in the right half-plane. By \eqref{eq:est-finite-comparison} and the half-plane Harnack inequality \eqref{equ:harmonic-plane}, we have
\[
G_{K_n}(0)=G_{K_n}\circ\mathcal{H}^{-1}(1)\ge\frac{G_{K_n}(x_0)}{\mathcal H(x_0)}
\succeq\frac1{1+\log(1/\kappa_n)}.
\]
By \eqref{eq:est-finite-target} and note that $\kappa_n\to0^+$, we have
\begin{equation}
G(t_n)=G_{\overline{D}_n}(0)-G_{K_n}(0)<0
\end{equation}
for all large $n$, which contradicts to the assumption that $G(t_n)\ge0$ for all $n$. Hence $K_n\subset W_0\cup\{-1\}$ for all large $n$. Without loss of generality, we assume that this holds for all $n$.

\medskip
We therefore have $K_n=\overline{h_0(D_n)}$, with the limit $K_*:=\overline{h_0(D_*)}$ in the sense of Carath\'{e}odory kernal convergence, where $D_*=\D(C_*,r_*)$. The circle $\partial D_*$ passes through $v$ and $0$ and meets $(-\infty,v]$ only at $v$. Since $g\circ h_0=\mathrm{id}$, the map $h_0$ is holomorphic near $0$ and has a square-root expansion at $v$, the circle parametrizations, composed with $h_0$, converge uniformly.
Thus there exists a constant $M'>1$ such that $\partial K_n$ are $M'$-quasi-circles. By \cite[Theorem 9.11, p.\,283]{Pom75},
\begin{equation}
\psi_{\overline{D}_n}^{-1}\to \psi_{\overline{D}_*}^{-1} \text{\quad and\quad} \psi_{K_n}^{-1}\to\psi_{K_*}^{-1}
\end{equation}
uniformly on $\EC\setminus\D$ as $t_n\to t_*$, where $\psi_{\overline{D}_n}$, $\psi_{\overline{D}_*}$ and $\psi_{K_n}$, $\psi_{K_*}$ are normalized Riemann maps (see \eqref{equ:psi-K}).
In particular,
\begin{equation}\label{eq:est-moving-point}
\psi_{\overline{D}_n}(0)\to\psi_{\overline{D}_*}(0)\in\T \text{\quad and\quad}
\psi_{K_n}(0)\to\psi_{K_*}(0)\in\T.
\end{equation}
Moreover,
\begin{equation}\label{equ:s-star}
e^{2\pi i s_n}:=\frac{\psi_{K_n}(0)}{|\psi_{K_n}(0)|}/\frac{\psi_{\overline{D}_n}(0)}{|\psi_{\overline{D}_n}(0)|}\to \frac{\psi_{K_*}(0)}{\psi_{\overline{D}_*}(0)}=:e^{2\pi i s_*}.
\end{equation}
Consequently the circle homeomorphisms satisfy (see \eqref{equ:g-t})
\[
 g_n:=\psi_{\overline{D}_n}\circ g\circ \psi_{K_n}^{-1}\circ R_{s_n}
 \to
 \psi_{\overline{D}_*}\circ g\circ \psi_{K_*}^{-1}\circ R_{s_*}=:g_*
\]
uniformly on $\T$.
Since the inverse $h_0:\overline{D}_*\to K_*$ of $g$ restricts to an orientation-preserving boundary homeomorphism, it follows that $g_*$ is again an orientation-preserving circle homeomorphism. Denote $\zeta_*:=\psi_{\overline{D}_*}(0)\in\T$. Since $g(0)=0$, by \eqref{equ:s-star} we have
\begin{equation}
g_*(\zeta_*)=\psi_{\overline{D}_*}\circ g\circ \psi_{K_*}^{-1}\circ R_{s_*}(\psi_{\overline{D}_*}(0))=\psi_{\overline{D}_*}\circ g(0)=\psi_{\overline{D}_*}(0)=\zeta_*.
\end{equation}
Thus $\zeta_*\in\T$ is a fixed point of $g_*$. Hence $\rot(g_*)=0$, and the continuity of the rotation number gives $\rho(t_n)\to0$.
\end{proof}

\subsection{Proof of the Main Lemma}

In this subsection, we combine Propositions \ref{prop:small-circle}, \ref{prop:est-infinity} and \ref{prop:est-finite-boundary} to prove Lemma \ref{lem:main-2} and hence the {\hyperref[lem:main]{Main Lemma}}.

\begin{proof}[Proof of Lemma \ref{lem:main-2}]
Let $\theta\ne0$ be a given irrational number. We define
\[
\xi:=e^{2\pi i\theta} \text{\quad and\quad} Z(t):=e^{G(t)+2\pi i\rho(t)}.
\]
It suffices to find a parameter $t_0\in\Omega\setminus\{0\}$ for which $Z(t_0)=\xi$.
Let $0<r_0<|v|/4$ be the small number introduced in Proposition \ref{prop:small-circle}.
On the parameter circle $\{t\in\C:|t|=r_0\}$, we have $G(t)>0$ and the winding number
\begin{equation}\label{eq:innerwind}
\operatorname{wind}(Z|_{|t|=r_0},\xi)=1.
\end{equation}
By Proposition \ref{prop:est-infinity}, there is $R_0>1$ such that $G(t)<0$ whenever $|t|\ge R_0$.
Choose $M>R_0+1$ and a small open arc $I\subset\T$ containing $0$ but not $\theta$.
By the compactness and Proposition \ref{prop:est-finite-boundary}, there exists a sufficiently small $\varepsilon>0$, such that the Jordan domain
\[
\Omega_\varepsilon:=\{t\in\C:|\im t|<M \text{ and }\re(t^2)<\tfrac{|v|}{2}-\varepsilon\}
\]
satisfies $\overline{\D}_{r_0}\subset \Omega_\varepsilon\Subset \Omega$ and
\[
t\in\partial \Omega_\varepsilon\,\Longrightarrow\, G(t)<0\ \text{or}\ \rho(t)\in I.
\]

On the boundary of $\Omega_\varepsilon$, we consider the radial homotopy
\[
Z_\sigma(t):=e^{(1-\sigma)G(t)-\sigma+2\pi i\rho(t)}, \text{\quad where } 0\le\sigma\le1.
\]
If $G(t)<0$ for some $t\in\partial\Omega_\varepsilon$, then $Z_\sigma(t)<1$ for all $0\le\sigma\le1$. If $\rho(t)\in I$ for some $t\in\partial\Omega_\varepsilon$, then the above homotopy still avoids $\xi$.
Its final image $Z_1(t)$ lies on the circle $\T_{1/e}$, and therefore
\begin{equation}\label{eq:outerwind}
\operatorname{wind}(Z|_{\partial \Omega_\varepsilon},\xi)=0.
\end{equation}
Orient both parameter boundary components of the closed annulus $\overline{\Omega}_\varepsilon\setminus\D_{r_0}$ counterclockwise.
If $Z(t)=\xi$ had no solution in $\overline{\Omega}_\varepsilon\setminus\D_{r_0}$, then the images of its inner and outer boundaries in $\C\setminus\{\xi\}$ would have the same winding number, which contradicts to \eqref{eq:innerwind} and \eqref{eq:outerwind}.
Hence there exists $t_0$ such that
\begin{equation}
G(t_0)=0 \text{\quad and\quad}\rho(t_0)=\theta.
\end{equation}
Then $\phi_{t_0}(0)=0$ by \eqref{equ:phi-t-G-equi}, and $g_{t_0}:\T\to\T$ in \eqref{equ:g-t} is orientation-preservingly conjugate to $g\circ\phi_{t_0}|_{\Gamma_{t_0}}$, which has rotation number $\theta$.
This completes the proof of Lemma \ref{lem:main-2} and hence the \hyperref[lem:main]{Main Lemma}.
\end{proof}

\bibliographystyle{amsalpha}
\bibliography{E:/Latex-model/Ref1}

\end{document}